\documentclass[twoside]{amsart}
\usepackage{amsmath,amssymb}

\newtheorem{theorem}{Theorem}
\newtheorem{lemma}[theorem]{Lemma} 
 
\newtheorem{corollary}[theorem]{Corollary}

\begin{document}
\title[Almost disjoint subgroups]{Almost disjoint pure subgroups of the
Baer-Specker group}  
\author{Oren Kolman}
\address{King's College London\\
Strand\\
London WC2R 2LS, UK}
\thanks{We thank the referee for constructive comments.}
\email{okolman@member.ams.org}

\author{Saharon Shelah}
\thanks{This research was partially supported by the German-Israel
Foundation for Science; publication number 683.}
\address{Institute of Mathematics\\
Hebrew University\\
Jerusalem, Israel}
\email{shelah@math.huji.ac.il}

\begin{abstract}
We prove in ZFC that the Baer-Specker group ${\bf Z}^\omega$ has $2^{\aleph_1}$
non-free pure subgroups of cardinality $\aleph_1$ which are almost disjoint:
there is no non-free subgroup embeddable in any pair.
\end{abstract}
\thispagestyle{empty}

\maketitle

In this short paper we prove the following result.

\begin{theorem}
\label{1}There exists a family ${\bf
G}=\{G_\alpha:\alpha<2^{\aleph_1}\}$ of non-isomorphic
non-free pure subgroups of the Baer-Specker group ${\bf Z}^\omega$ such that:\hfil\break
(1.1)\quad each $G_\alpha$ has cardinality $\aleph_1$;\hfil\break
(1.2)\quad if $\alpha<\beta$, then $G_\alpha$ and $G_\beta$ are almost disjoint:
if $H$ is isomorphic to subgroups of $G_\alpha$ and $G_\beta$, then $H$ is free.
In particular, $G_\alpha\cap G_\beta$  is free.
\end{theorem}
Recall that the Baer-Specker group ${\bf Z}^\omega$ is the abelian group of
functions from the natural numbers into the integers (see \cite{1} and \cite{18}). It contains
the canonical pure free subgroup ${\bf Z}_\omega=\oplus_{n<\omega}{\bf Z}$.
The group ${\bf Z}^\omega$ is not $\kappa$-free for any cardinal $\kappa>\aleph_1$,
but it is $\aleph_1$-free, so the groups $G_\alpha$ in Theorem \ref{1} are almost free.

Theorem \ref{1} answers a question of the first author, and has its place in the
line of recent research dealing with the lattice structure of the pure subgroups
of ${\bf Z}^\omega$ (see \cite{2}, \cite{3}, and \cite{5}--\cite{8}). For example, Irwin
asked whether there is a subgroup of ${\bf Z}^\omega$ with uncountable dual but
no free summands of infinite rank. This problem was resolved recently by Corner
and Goebel \cite{5} who proved the following stronger fact.

\begin{theorem}\cite{5}
The Baer-Specker group ${\bf Z}^\omega$ contains a pure subgroup $G$ whose endomorphism
ring splits as End$(G)={\bf Z}\oplus$Fin$(G)$, with $|G^*|=2^{\aleph_0}$, where $
{\bf Z}$ is the scalar multiplication by integers and Fin$(G)$ is the ideal of
all endomorphisms of $G$ of finite rank.
\end{theorem}
Quotient-equivalent and almost disjoint abelian groups have been studied by
Eklof, Mekler and Shelah in \cite{9}--\cite{11}, who showed that under various
set-theoretic hypotheses, there exist families of maximal possible size of almost
free abelian groups which are pairwise almost disjoint. Following \cite{11}, we say
that two groups A and B are almost disjoint if whenever H is embeddable as a subgroup
in both A and B, then H is free. Clearly if A and B are non-free and almost disjoint,
then they are non-isomorphic in a very strong way. On the other hand, the intersection
of two almost disjoint groups of size $\aleph_1$ need not necessarily be countable, so
group-theoretic almost disjointness differs from its set-theoretic homonym. Theorem \ref{1}
establishes in ZFC that the Baer-Specker group contains large families of almost disjoint
almost free non-free pure uncountable subgroups.

Our group and set-theoretic notation is standard and can be found in \cite{10} and
\cite{14}. For example, ${}^{\omega_1>}2$ is the set of partial functions from
$\omega_1$ into $\{0,1\}$ whose domains are at most countable; ${}^{\omega_1}2$
is the set of all functions from $\omega_1$ into $\{0,1\}$; for a regular cardinal
$\chi,\ H(\chi)$ is the family of all sets of hereditary cardinality less than
$\chi$.

For a set $A\subseteq H(\chi)$ for $\chi$ large enough, we write dcl$_{\big(H(
\chi),\in,<\big)}[A]$ for the Skolem closure (Skolem hull) of $A$ in
the structure $\big(H(\chi),\in,<\big)$, where $<$ is a well-ordering
of $H(\chi)$ (for details, see \cite{16}, 400-402, or \cite{15}, 165-170).

In proving Theorem \ref{1}, we shall appeal to the well-known Engelking-Kar\l owicz
theorem from set-theoretic topology:

\begin{theorem} \cite{13}
If $|Y|=\mu=\mu^{<\sigma}<\lambda=|X|\leq2^\mu$, then there are functions $h_
\alpha:X\rightarrow Y$ for $\alpha<\mu$ such that for every partial function
$f$ from $X$ to $Y$ of cardinality less than $\sigma$, for some $\alpha<\mu,
\ f\subseteq h_\alpha$.
\end{theorem}
A self-contained short proof can be found in \cite{17}, 422-423. We shall need
just the case when $\mu=\sigma=\aleph_0$, and $\lambda=2^\mu$. Since it may be
less familiar to algebraists, for convenience we deduce the fact to which we appeal
later on (although it also appears as Corollary 3.17 in \cite{4}).

\begin{lemma}
\label{2}There exists a family $\{f_\eta:\eta\in^{\omega_1>}\!2\}$ such that $f_\eta:\omega
\rightarrow{\bf Z}$, and whenever $\eta_1,\dots,\eta_k$ are distinct and $a_1,
\dots,a_k\in{\bf Z}$, then\hfil\break
$\{i<\omega:(\forall\ l\leq k)(f_{\eta_l}(i)=a_l)\}$ is infinite.
\end{lemma}

\begin{proof} Take $\mu=\sigma=\aleph_0,\ \lambda=2^\mu,\ X={}^{\omega_1>}2$ and $Y={\bf Z}$
in the Engelking-Kar\l owicz theorem. Since $|^{\omega_1>}2|=2^{\aleph_0}$ and $|{\bf
Z}|=\aleph_0$, we know that there exist functions $h_n:{}^{\omega_1>}2\rightarrow{\bf
Z}$ for $n<\omega$ such that for every partial function $f$ from ${}^{\omega_1>}2$ to
${\bf Z}$ whose domain is finite, there is some $m<\omega$ such that $f\subseteq h_m$.
Let $\{g_i:i<\omega\}$ be an enumeration with infinitely many repetitions of each $h_n$
for $n<\omega$.

For each $\eta\in^{\omega_1>}2$, define $f_\eta:\omega\rightarrow{\bf Z}$ by
$f_\eta(i)=g_i(\eta)$. The family $\{f_\eta:\eta\in^{\omega_1>}2\}$ is as required:
for if $\eta_1,\dots,\eta_k$ are distinct and $a_1,\dots,a_k\in{\bf Z}$ are
given, then the set $f=\langle(\eta_1,a_1),\dots,(\eta_k,a_k)\rangle$
is a finite function,
so there is some $m$ such that $f\subseteq h_m$ and it is now easy to see that
$\{i<\omega:(\forall\ l\leq k)(f_{\eta_l}(i)=a_l)\}$ is infinite.
\end{proof}

A well-known algebraic fact will also be useful:

\begin{lemma}
\label{3}Let $C$ be a closed unbounded subset of the regular uncountable cardinal $\kappa$.
Suppose that $H$ is an abelian group of cardinality $\kappa$, and
$\langle H_\alpha:
\alpha<\kappa\rangle$ is a $\kappa$-filtration of $H$ (a continuous
increasing chain of subgroups
$H_\alpha,\ |H_\alpha|<\kappa$, whose union is $H$). Let $S=\{\alpha\in C:H/H_\alpha$ is not $\kappa$-free$\}$.
Then $H$ is free if and only if $S$ is non-stationary in $\kappa$.
\end{lemma}

\begin{proof} Well-known: see Proposition IV.1.7 in \cite{10}.
\end{proof}

We refer the reader to \cite{14} for the definitions of the
characteristic $\chi(g)$ and the type $\tau(g)$
of an element $g$ in a group.

\vskip 10pt
\noindent Now we prove Theorem \ref{1}.

\begin{proof}
Let ${\bf P}$ be the set of prime numbers, and let $\{P_\eta:\eta\in^{\omega_
1>}2\}$ be a family of almost disjoint (infinite) subsets of $\bf P$:
$\eta\ne\nu\in^{\omega_1>}2\ \Rightarrow\ |P_\eta\cap P_\nu|<\aleph_0$.
By Lemma \ref{2}, there exists $\{f_\eta:\eta\in^{\omega_1>}2\}$ such that
$f_\eta:\omega\rightarrow{\bf Z}$, and if $\eta_1,\dots,\eta_k$ are distinct and $a_1,\dots,a_k\in{\bf Z}$,
then $\{i<\omega:(\forall\ l\leq k)(f_{\eta_l}(i)=a_l)\}$ is infinite.

Define functions $x_\eta$ and $x_{\eta,j}$ in ${\bf Z}^\omega$ as follows.
Let $x_\eta=\langle\pi_{\eta,i}\cdot f_\eta(i):i<\omega\rangle$ where
$\pi_{\eta,i}=\Pi\{p\in P_\eta:p<i\}$,
and let $x_{\eta,j}=\langle\pi^j_{\eta,i}\cdot
f_\eta(i):i<\omega\rangle$ where $\pi^j_{\eta,i}
=\Pi\{p\in P_\eta:j\leq p<i\}$ (=0 if $i\leq j$). Note that
$x_\eta=x_{\eta,0}$.

For $\eta\in^{\omega_1}2$, let $G_\eta$ be the subgroup of ${\bf Z}^\omega$
generated by ${\bf Z}_\omega\cup\{x_{\eta|\alpha,j}:\alpha<\omega_1,0\leq j<\omega\}$.

We show that the family ${\bf G}=\{G_\eta:\eta\in^{\omega_1}2\}$ satisfies
the conclusions of Theorem \ref{1}.
\vskip 10pt
\noindent{\bf Claim 1:} $G_\eta$ is pure in ${\bf Z}^\omega$.
\vskip 10pt
\noindent{\bf Proof of Claim 1:} Suppose that $rx=g$ for some $x\in{\bf Z}^\omega,\ r\in{\bf
N}$, and $g\in G_\eta$. Say $g=y+n_1x_{\eta|\alpha_1,j_1}+\dots+n_mx_{\eta|\alpha
_m,j_m},\ n_l\ne0$, with $y\in{\bf Z}_\omega$. Without loss of generality (adding
more elements from ${\bf Z}_\omega$ to the RHS if necessary), $(\forall\ l\leq
m)(j_l=j)$ for some $j<\omega,\ j>r,\ y(i)=0\ (\forall\ i>j)$, and $x(i)=0\ (\forall
i\leq j)$. Relabelling (if necessary), we may assume that $\alpha_1<\dots<\alpha
_m<\omega_1$, and because $x_{\eta|\alpha_l,j}(i)=0$ if $i\leq j$, we may write
$$
rx=ry^*+n_1x_{\eta|\alpha_1,j}+\dots+n_mx_{\eta|\alpha_m,j},\quad{\rm for\ some}
\ y^*\in{\bf Z}_\omega.
$$
Fix $k\in\{1,\dots,m\}$. Since $\eta|\alpha_1,\dots,\eta|\alpha_m$ are distinct
$(\alpha_1<\dots<\alpha_m)$, letting $a_l=\delta_{kl}$ (Kronecker delta), we
know that the set $N_k=\{i<\omega:(\forall\ l\ne k)(f_{\eta_l}(i)=0,\ f_{\eta_k}(i)
=1)\}$ is infinite. For large enough $i$ in this set \big(e.g. $i>{\rm max}_{1\leq l
\leq m}\ [{\rm min}(P_{\eta|\alpha_l}\backslash\{0,\dots,j\})]\big),\ x_{\eta|
\alpha_l,j}(i)$ is zero if and only if $l\ne k$. So for infinitely many $i<\omega$,
for $l\ne k,\ x_{\eta|\alpha_l,j}(i)=0$, and $x_{\eta|\alpha_k,j}(i)\ne0$.

Unfix $k$. For each $k\leq m$, for infinitely many $i\in(j,\omega)\cap N_k,\ rx(i)=n_k
x_{\eta|\alpha_k,j}(i)=n_k\Pi\{p\in P_{\eta|\alpha_k}:j\leq p<i\}$. Since $r<j$, we must have $rs_k=
n_k$ for some $s_k$ in $\bf Z$, and therefore $x=y^*+s_1x_{\eta|\alpha_1,j}+\dots
+s_mx_{\eta|\alpha_m,j}\in G_\eta\ (G_\eta$ is torsion-free). Hence $G_\eta$ is
pure in ${\bf Z}^\omega$, which establishes Claim 1.
\vskip 10pt
\noindent{\bf Claim 2:} $G_\eta$ has cardinality $\aleph_1$, so (1.1) holds.
\vskip10pt
\noindent{\bf Proof of Claim 2:} If $\xi\ne\zeta\in^{\omega_1>}2$, then for some $j<\omega,\ P_
\xi\cap P_\zeta\subseteq j$. Pick $p,q>j$ with $p\in P_\xi$ and $q\in P_\zeta$;
so the set $B=\{i<\omega:f_\xi(i)=p$ and $f_\zeta(i)=q\}$ is infinite, and if
$i\in B$ is bigger than max$\{j,p,q\}$, then $x_{\xi,j}(i)\ne x_{\zeta,j}(i)$,
since $x_{\xi,j}(i)$ is non-zero and divisible by $p^2$ but by no prime in $P_
\zeta$, and $x_{\zeta,j}(i)$ is non-zero and divisible by $q^2$ but by no prime
in $P_\xi$. It follows that $G_\eta$ has cardinality $\aleph_1$.
After this observation, a second's reflection on the element types of
$G_{\eta}$ and $G_{\nu}$ (for $\eta\ne\nu$)
should convince the reader that the groups are neither isomorphic nor free.
\vskip 10pt
\noindent{\bf Claim 3:} (1.2) holds: if $\eta_1\ne\eta_2\in^{\omega_1}2$, then $G_{\eta_1}$
and $G_{\eta_2}$ are almost disjoint.
\vskip 10pt
\noindent{\bf Proof of Claim 3:} Suppose (towards a contradiction) that for some $\eta_1\ne\eta
_2\in^{\omega_1}2$, for some non-free abelian group $H$, there exist isomorphisms
$\varphi_l:H\rightarrow$ range$(\varphi_l)\leq G_{\eta_l},\ l=1,2$. Since $G_{\eta_l}$
is $\aleph_1$-free, $H$ must have cardinality $\aleph_1$. Let
$\langle H_i:i<\omega_1\rangle$ be an
$\omega_1$-filtration of $H$. Without loss of generality, we may assume that each
$H_i$ is pure in $H$, so that $H/H_i$ is torsion-free.

Let $G_{\eta,i}=\langle{\bf Z}_\omega\cup\{x_{\eta|\beta,j}:j<\omega,\ \beta<i\}
\rangle$ for $i<\omega_1$ and $\eta\in\{\eta_1,\eta_2\}$.

Note that $\langle G_{\eta,i}:i<\omega_1\rangle$ is a
$\omega_1$-filtration of $G_\eta$,
since it is increasing and continuous with union $G_\eta$, and each
$G_{\eta,i}$ is countable. For large enough $\chi$,
the set $C$ defined by $\{\delta<\omega_1:{\rm dcl}_{\big(H(\chi),\in,<\big)}
[\delta\cup\{G_{\eta_1},G_{\eta_2},\{x_\nu,f_\nu:\nu\in^{\omega_1>}2\},\eta_1,
\eta_2,\varphi_1,\varphi_2,\{H_i:i<\omega_1\}\}]\cap\omega_1=\delta\}$ is a club of $\omega_1$
(well-known, or see \cite{16}, 401). Note that if $\delta\in C$, then $\varphi_l$ maps $H_\delta$
into $G_{\eta_l,\delta}$. Since $H$ is not free, it follows by Lemma \ref{3} that $S=
\{\delta\in C:H/H_\delta$ is not $\aleph_1$-free$\}$ is stationary. By Pontryagin's Criterion,
for each $\delta\in S,\ H/H_\delta$ has a non-free (torsion-free) subgroup $K_
\delta/H_\delta$ of finite rank $n_\delta+1$ such that every subgroup of $K_\delta
/H_\delta$ of rank less than $n_\delta+1$ is free. Let $H_\delta^{\ +}/H_\delta$
be a pure subgroup of $K_\delta/H_\delta$ of rank $n_\delta$. Then $H_\delta^{\ +}
/H_\delta$ is free with basis $y_0+H_\delta,\dots,y_{n_\delta-1}+H_\delta$ say.
So $K_\delta/H_\delta^{\ +}\simeq(K_\delta/H_\delta)/(H_\delta^{\ +}/H_\delta)$
is a torsion-free rank-1 group which is not free, and hence there is a non-zero
element $y_{n_\delta}+H_\delta^{\ +}$ which is divisible in $K_\delta/H_\delta^
{\ +}$ by infinitely many natural numbers. Call this set of natural numbers $A$.

For $l=1,2$, for large enough $j_l(*)<\omega$, and $\beta^l_{\ 0}<\dots<\beta^l
_{\ k_l}<\omega_1,\ \varphi_l(y_m)$ is an element of the subgroup of $G_{\eta_l}$
generated by $G_{\eta_l,\delta}\cup\{x_{\eta_l|\beta^l_{\ 0},j_l(*)},\dots,x_{
\eta_l|\beta^l_{\ k_l},j_l(*)}\}$ for all $m\leq n_{\delta}$.

Taking large enough $\delta\in S$, we may assume that min$\{\alpha:\eta_1|\alpha
\ne\eta_2|\alpha\}<\beta^l_{\ 0},\ l=1,2$. Since $\delta\in C$, we can show the
following claims:\hfil\break
$(*)_1$: The set $A$ does not contain infinitely many powers of one prime.\hfil\break
$(*)_2$: The set $Q=({\bf P}\cap A)\subseteq P_{\eta_l|\beta^l_{\ 0}}\cup\dots\cup P
_{\eta_l|\beta^l_{\ k_l}}$.

Now $(*)_1$ is true because non-zero sums of elements in\hfil\break
$G_{\eta_l,\delta}\cup\{x_{\eta_l|\beta^l_{\ 0},j_l(*)},\dots,x_{\eta_l|\beta^l_{\ k_l},j_l(*)}\}$ are divisible
by at most finitely many powers of any given prime (by the definition of the elements
$x_{\eta_l|\beta,j})$. Note that $\chi(y_{n_\delta}+H_\delta^{\ +})=\cup_{\{y\in y_{n_
\delta}+H_\delta^{\ +}\}}\chi(y)\leq cup_{\{y\in y_{n_
\delta}+H_\delta^{\ +}\}}\chi(\varphi_l(y))$, where the characteristics are taken relative to
$K_\delta/H_\delta^{\ +}$, $K_\delta$ and\hfil\break
$G_{\eta_l,\delta}\cup\{x_{\eta_l|\beta^l_{\ 0},j_l(*)},\dots,x_{\eta_l|\beta^l_{\ k_l},j_l(*)}\}$ respectively.
Hence $(*)_1$ holds. By $(*)_1$, since $A$ is infinite, the set $Q={\bf P}\cap A$ is infinite.

Also, the same characteristic inequality implies that $Q\subseteq P_{\eta_l|\beta^l_{\ 0}}\cup\dots\cup P_{\eta_l|\beta^l_{\ k
_l}}$. So $(*)_2$ is true. Hence, $Q\subseteq\cap_{l=1,2}(\cup_{k\leq k_l}P_{\eta_l|\beta^l_{\ k}})$
which is finite (since the family $\{P_\eta:\eta\in^{\omega_1>}2\}$
is almost disjoint).
This is a contradiction, and so Claim 3 follows, completing the proof
of Theorem \ref{1}.
\end{proof}

\begin{corollary}
Every non-slender $\aleph_1$-free abelian group $G$ has a family $\{G_\alpha:\alpha<2^{\aleph_1}\}$
of non-free subgroups such that:\hfil\break
1. each $G_\alpha$ is almost free of cardinality $\aleph_1$;\hfil\break
2. if $\alpha<\beta$, then $G_\alpha$ and $G_\beta$ are almost disjoint.
\end{corollary}
\begin{proof} By Nunke's characterisation of slender groups (see Corollary IX.2.5 in
\cite{10} for example), $G$ must contain a copy of the Baer-Specker group.
\end{proof}
\vskip10pt
\noindent{\bf Remark}: For the same reason, the corollary is true for any non-slender cotorsion-free
abelian group.

\end{document}